\documentclass[final]{siamltex}

\title{Obstructions to deformation of curves to other hypersurfaces}

\author{ B. Wang}

\begin{document}

\maketitle

\begin{abstract}
We are interested in obstructions to the FIRST order deformation of a pair of a  smooth hypersurface $f_0$ and a smooth curve $C_0$  contained in $f_0$.  
In the first half of the paper, we
give necessary conditions for the pair to deform in the first order.
In particular, for a rational curve $C_0$, this necessary condition is
$$H^1(N_{C_0}f_0(1))=0.$$  In the second half, we apply the necessary conditions from the first half of the paper  to study the
geometry of  smooth curves in hypersurfaces (theorem 4.1,  theorem 5.1).  
The main application is for the case where $C_0$ is a rational curve.

\end{abstract}

\pagestyle{myheadings}
\thispagestyle{plain}
\markboth{ B. Wang}{An obstruction to the deformations of curves to all hypersurfaces}

\section{Introduction} Let $\mathbf P^n$ be the projective space of dimension $n\geq 3$ over complex numbers. 
Let $f_0\subset \mathbf P^n$ be a smooth hypersurface of degree $h$. 
 Let $C_0\subset f_0$ be a smooth curve. We investigate the existence of a family of pairs
$C_t\subset  f_t$, the curves $C_t$ of degree $d$ and the hypersurfaces $f_t$ of degree $h$ in the projective space
$\mathbf P^n$, where $t$ is in a variety.
 Furthermore the family $f_t$ is not a constant hypersurface. We simply call this $C_t$ 
a deformation of $C_0$ to other hypersurfaces, or a ``FULL" deformation of the pair to other hypersurface. The similar question was investigated by L. Chiantini and Z. Ran 
in [5]. In general there is  Kodaira's deformation theory ([6], theorem 1)
about the submanifold $C_0\subset f_0$ in a complex manifold $f_0$, that says a sufficient condition for the 
$C_0$ to deform to ALL the other submanifolds  is $$H^1(N_{C_0}f_0)=0.$$ The space $H^1(N_{C_0}f_0)$ is often called
the obstruction space to the deformation of the pair $C_0\subset f_0$. 
But in general, it is not clear that this condition is also a necessary condition, i.e,  if $$H^1(N_{C_0}f_0)\neq 0,$$ 
$C_0$ may still be able to deform to all the other hypersurfaces (We don't have a proof of that yet). 
  So in the first half of the paper, we prove theorem  1.1  below that gives  necessary conditions (i.e. obstructions) for $C_0$ to deform to ``other hypersurfaces" in the first order. The conditions, the formula (1.6) or (1.10),  are expressed in terms of the dimensions of cohomology groups of various twisted bundles over $C_0$. But our assumption in theorem  1.1
is weaker than the existence of the deformation $C_t\subset f_t$.  We only use 
the first order deformation of $C_0\subset f_0$.\par
The second half of the paper concentrates on the applications. They can be categorized into two different kinds: \par
(1) Rational curves in a smooth quintic 3-fold.  In this case,  we see
the converse of the Kodaira's theorem ([6], theorem 1) for a rational curve $C_0$  in a generic quintic 3-fold $f_0$  is exactly 
the Clemens' conjecture: the $H^1$ of the normal bundle of the rational curve $C_0\subset f_0$ is equal to $0$, 
$$ H^1(N_{C_0}f_0)=0.$$
Thus even though our theorem  1.1 did not completely solve the Clemens' conjecture,  it still gives a new approach and a new result towards the Clemens' conjecture.  See theorem 4.1.\par
(2) Numerical bounds for curves in hypersurfaces. In general there is the interest in knowing
the lower bound of the geometric genus of a subvariety in a general hypersurface or a complete intersection, and numerical bounds of
other invariants. 
There are many works on this by Clemens, Chiantini, Ein, Pacienza, Ran, Voisin and Xu, etc. Among them some bounds are sharp. 
Even though this paper did not produce better bounds, but our results are new in the sense that we assumed a weaker condition: our
$f_0$ is not generic and the deformation of the pair of $C_0, f_0$ is in the first order only, i.e. we did not assume the pair
$C_0\subset f_0$ can actually deform to generic hypersurfaces. Please see the details in theorem 5.1, theorem 5.2. The difference between our bounds and their sharp bounds may be caused by our weaker assumption for the deformation of the pair (in the first order only ).  \bigskip

\bigskip

{\bf The formal setting}.\par
To state the theorem in a precise way, we need to give a formal description about the first order assumption.
Let $H^0(\mathcal O_{\mathbf P^n}(h))$ denote the vector space of homogeneous polynomials of degree $h$ in $n+1$ variables.
We use the same letter $f_0\in H^0(\mathcal O_{\mathbf P^n}(h))$ to denote the hypersurface $div(f_0)\subset \mathbf P^n$,
homogeneous polynomial $f_0$, and its projectivization in  $\mathbf P(H^0(\mathcal O_{\mathbf P^n}(h)))$. 
Let $S\subset \mathbf P(H^0(\mathcal O_{\mathbf P^n}(h)))$ be a subvariety containing $f_0$ which is
a smooth point of $S$. Also assume that $f_0$ is a smooth hypersurface.  Let 
\begin{eqnarray} && X_S\subset S\times \mathbf P^n,\\&&
X_S=\{(f, x): f\in S, f(x)=0\}.\end{eqnarray}
be the universal hypersurface.
\par
Let $C$ be a smooth projective curve of genus $g$, and
$$c_0: C\to   f_0\subset P^n$$
a smooth imbedding of $C$ to $f_0$. Then 
$$\bar c_0: C\to \{f_0\}\times f_0\subset X_S$$ is the induced imbedding.
The projection $$P_S: X_S\to S$$ induces a map on the sections of bundles over $C$, 
\begin{equation} P_S^s: H^0(\bar c_0^\ast(TX_S)) \to T_{f_0}S,\end{equation}
where $ T_{[f_0]}S\simeq H^0(T_{[f_0]}S\otimes \mathcal O_C)$ is the space of global sections of the 
 trivial bundle whose each fibre is  $T_{f_0}S$.
\medskip
In this paper we consider two specific parameter spaces for $S$:\smallskip

{\bf Assumption (1) } The first subvariety $S$ under consideration  is the collection of hypersurfaces 
in the following form:
\begin{equation} f_0+\sum_{i=0}^h a_i L_0\cdots\hat L_i\cdots L_{h},\quad (\hat L_i \ is \ omitted)\end{equation}
where $L_i\in H^0(\mathcal O_{\mathbf P^n}(1)), i=0, \cdots, h$ are fixed sections whose zeros are distinct, i.e.
$$div(L_i)\neq div(L_j), i\neq j. $$
Let $$A'=\mathbf C^{h+1}=\{(a_0, \cdots, a_h)\}$$ be the parameter space of the family. Let
$A\subset A'$ that parametrizes smooth hypersurfaces. So $S=A$ in this case.\medskip

{\bf Assumption (2) }  Secondly $S$ is the entire space $\mathbf P(H^0(\mathcal O_{\mathbf P^n}(h)))$. We will denote
$\mathbf P(H^0(\mathcal O_{\mathbf P^n}(h)))$ by $E$.  So $S=E$ in this case.
\bigskip

Let $N_{c_0}V$ denote the pullback of a bundle  $V\subset T\mathbf P^n|_{C_0}$ over $C_0$. 
Let  $h^i(E)$ denote the dimension of $H^i(E)$ for any sheaf $E$.
Let $\mathcal L=c_0^\ast(\mathcal O_{\mathbf P^n}(1))$. Then $d=deg(\mathcal L)$, $g=genus(C)$.

\bigskip

 \begin{theorem}   Assume $P_A^s$ is  surjective (see assumption (1) and formula (1.3)). 
\par
(1) If $C$ is a rational curve and
\begin{equation} \{L_i=0\}\cap \{L_i=0\}\cap C_0=\emptyset, i\neq j, \end{equation}
 then
 
\begin{equation} H^1(N_{c_0}f_0(1))=H^1(c_0^\ast (Tf_0(1)))=0,\end{equation}
where $N_{c_0}f_0(1)$ is the pull-back of the twisted normal bundle $N_{C_0}f_0(1)$. \par

(2) Let  \begin{eqnarray} \sigma(c_0, f_0)=&&(h+1)h^0(\mathcal L)-h-h^0(\mathcal L^{h+1})\\ &&
+h^1(c_0^\ast(Tf_0(1))) -h^1( c_0^\ast(T\mathbf P^n(1))+h^1(\mathcal L^{h+1}).\end{eqnarray}

If $$(h+1)h^0(\mathcal L)-h\geq h^0( c_0^\ast(T\mathbf P^n(1)),$$
and \begin{equation} \{L_i=0\}\cap \{L_i=0\}\cap C_0=\emptyset, i\neq j, \end{equation}
then
\begin{equation} \sigma(c_0, f_0)=0.\end{equation}

\end{theorem}

\bigskip

{\bf Remark} This theorem proves that if $c_0, f_0$ satisfy  conditions 
$$(h+1)h^0(\mathcal L)-h\geq h^0( c_0^\ast(T\mathbf P^n(1)), \quad and \ \sigma(c_0, f_0)\neq 0$$  
then $C_0$ can't deform to all the hypersurfaces in $A$ in 
the first order.  Thus
$H^1(N_{C_0}f_0(1))\neq 0$ in the case $g=0$, and $\sigma(c_0, f_0)\neq 0$ in general, give us  obstructions to the 
deformation of $C_0$ to other hypersurfaces. Note $ \sigma(c_0, f_0)=H^1(N_{C_0}f_0(1))$ in the case $g=0$. 

\bigskip

The rest of the paper is organized as follows. In section 2, we study the deformation of a smooth hypersurface. We mainly describe and prove a theorem by H. Clemens, then some sequences of bundles associated with it. This is the main technique for
the paper. In section 3, we study the first order deformation of smooth curves together with hypersurfaces containing it. This is a proof of theorem  1.1. In 
section 4, we apply theorem 1.1 to rational curves in quintic 3-folds. The main result is to determine the upper-bound of the degrees of the summand in the normal bundle
of a rational curve in  a smooth quintic 3-fold. In section 5, we apply theorem  1.1 to any smooth hypersurfaces. We obtain bounds for a couple of numerical invariants. These
bounds are not better than the bounds in Clemens' paper [4], however they are obtained by using our weaker first order assumption. Therefore they are new theorems (theorem
5.1, theorem 5.2) .

\section{Deformation of hypersurfaces} The main idea of the proof is to transform the problems
of $T\mathbf P^n$ to similar types of problems of some isomorphic bundle ${TX_A(1)\over G(1)}$. Then the existence of the
first order deformation of the pair $C_0, f_0$ allows us to work with 
${TX_A(1)\over G(1)}$, which is more accessible now than $T\mathbf P^n$. 
 Thus the isomorphism between 
$T\mathbf P^n$ and ${TX_A(1)\over G(1)}$ serves as an important bridge between two different realms. 
In this section, we introduce the construction of the vector bundle ${TX_A(1)\over G(1)}$ and the 
associated morphisms that are  used in our proof. \par
Let \begin{equation} F(a_1, \cdots, a_h, x)=f_0(x)+\sum_{i=0}^h a_i L_0(x)\cdots\hat L_i(x)\cdots L_{h}(x), \quad (omit \ L_i)\end{equation}
be the universal polynomial. Thus
$$\{F=0\}=X_A\subset A\times \mathbf P^n.$$
is also the universal hypersurface, which is smooth. 
Let $W\subset \mathbf P^n$ denote the complement of the proper subvariety
$$\cup_{h\geq  j>i\geq 0}\{L_i=L_j=0\}.$$
Let \begin{eqnarray} && X_W=X_A\cap (A\times W)\\&&
f_0^W=f_0\cap W.\end{eqnarray}
Let
\begin{equation} u_i=L_0{\partial\over \partial a_0}-L_i{\partial\over \partial a_i}, i=1, \cdots, h
\end{equation} be sections
of $ TA\otimes \mathcal O_{W}(1)$. Since $u_i$ annihilate $F$, they are tangent to
$X_W$. So let 
\begin{equation} G(1)\subset  TX_W(1)\end{equation}
be the vector bundle of rank $h$ over $X_W$ that is generated by the sections
$u_i$. 
 
\bigskip
For any smooth varieties $V_1, V_2$, let
$$T_{V_1/V_2}$$ denote the relative tangent bundle of $V_1$ over $V_2$, i.e. it is
the bundle $TV_1\oplus \{0\}$ over the variety $V_1\times V_2$. 
\bigskip

The following theorem 2.1  is communicated to us by H. Clemens([3]). 

\bigskip

\begin{theorem} ({\rm H. Clemens})
\begin{equation} {TX_W(1)\over G(1)}\simeq T_{W/A}(1),
\end{equation}
where $T_{W/A}(1)$ is restricted to $X_W$. 

\end{theorem}
\begin{proof} Consider the exact sequence
\begin{equation}\begin{array}{ccccccccc}
0&\rightarrow {TX_W(1)\over G(1)}&\rightarrow {T(A\times W)(1)\over G(1)} 
&\rightarrow &\mathcal D &\rightarrow 0.\end{array}\end{equation}
of bundles over $X_W$, 
where $\mathcal D$ is some quotient bundle over $X_W$.
Easy to see \begin{equation} c_1(\mathcal D)=c_1(\mathcal O_{\mathbf P^n}(h+1))|_{X_w}.\end{equation}
Let $s$ be a generic section of $\mathcal O_{\mathbf P^n}(1)$ that does not have common zeros with
$L_i, i=0, \cdots, h$. Let $\sigma$ be the reduction of  $s{\partial \over \partial a_0}$
in ${T(A\times W)(1)\over G(1)}$. Notice the zeros of $\sigma$ is exactly
\begin{equation} div(\sigma)=div(s L_1\cdots L_h).\end{equation}
Since $s L_1\cdots L_h\in H^0(\mathcal O_{\mathbf P^n}(h+1))$, $\sigma$ splits the sequence (2.7).
If $L_s\subset {T(A\times W)(1)\over G(1)}$ is the line bundle generated by $\sigma$, 
\begin{equation} L_s\oplus {TX_W(1)\over G(1)}={T(A\times W)(1)\over G(1)},\end{equation}
as bundles over $X_W$.
Secondly, we have another exact sequence
\begin{equation}\begin{array}{ccccccccc}
0&\rightarrow T_{W/A}(1)&\rightarrow {T(A\times W)(1)\over G(1)} 
&\rightarrow &\mathcal D' &\rightarrow 0.\end{array}\end{equation}
of bundles over $X_W$,
where $\mathcal D'$ is some quotient bundle over $X_W$.
By the direct calculation (note $G(1)$ is a trivial bundle):
$$c_1(\mathcal D')=c_1(c_0^\ast(T_{A/W}(1)))=(h+1)(c_1(\mathcal O_{\mathbf P^n}(1)))|_{X_W}$$
As above,  $\sigma$ splits this sequence (2.11). Hence
\begin{equation} L_s\oplus T_{W/A}(1)={T(A\times W)(1)\over G(1)}.\end{equation}
Comparing (2.10), (2.12), we obtain
\begin{equation} {TX_W(1)\over G(1)}\simeq T_{ W/A}(1), \end{equation}
over $X_W$.\quad 
\end{proof}

\bigskip

 Let 
\begin{equation}\begin{array}{ccc} I: {TX_W(1)\over G(1)} &\rightarrow & T_{ W/A}(1) \\
({\partial \over \partial a_0}, v) &\rightarrow  & -v\end{array}\end{equation}
be this isomorphism in the formula (2.13) or (2.6),
where $({\partial \over \partial a_0}, v)$ is the decomposition in the formula (2.12) and $v\in T_xW$.

Consider the composition map $\mu_1$: 
\begin{equation}\begin{array}{ccccc} T_{A/W} &\rightarrow & {T(A\times W)\over TX_W} 
&\rightarrow  & N_{f_0^W} W\simeq \mathcal O_{\mathbf P^n}(h)|_{f_0^W}\\
\alpha &&\rightarrow && {\partial F\over \partial \alpha}|_{\bar f_0^W}
\end{array}\end{equation}
where $\bar f_0^W=\{f_0\}\times f_0^W$. 
The last map is the restriction map. Tensoring it  with  $\mathcal O_{\mathbf P^n}(1)$, we have the composition $\mu_2$:
\begin{equation}\begin{array}{ccccccc} TX_W(1) &\rightarrow  T_{A/W}(1) & \rightarrow &
\mathcal O_{\mathbf P^n}(h+1)|_{f_0^W}
\end{array}\end{equation}
where the first map is the differential of the projection $X_W\to A$. 
Since $\mu_2$ vanishes on $G(1)$, we obtain the bundle morphism $\mu_3$:
\begin{equation}\begin{array}{ccc} {TX_W(1)\over G(1)}|_{\bar f_0^W}
 &\stackrel{\mu_3} \rightarrow  \mathcal O_{\mathbf P^n}(h+1)|_{f_0^W}.
\end{array}\end{equation}

\begin{lemma} There is a commutative diagram
\begin{equation}\begin{array}{ccccccccc}
0 & \rightarrow &T_{f_0^W/A}(1) &\rightarrow & {TX_W(1)\over G(1)}|_{\bar f_0^W}
&\stackrel{\mu_3}\rightarrow & \mathcal O_{\mathbf P^n}(h+1)|_{f_0^W}& \rightarrow & 0\\
& & \downarrow\scriptstyle{I} && \downarrow\scriptstyle{I} && \|&&\\
0 & \rightarrow &Tf_0^W (1)&\rightarrow & {TW(1)}|_{f_0^W}
& \stackrel{\nu}\rightarrow &\mathcal O_{\mathbf P^n}(h+1)|_{f_0^W}& \rightarrow &0
\end{array}\end{equation}
where $I$ is the isomorphism in the formula (2.14) and
$\nu$ is the differential map
$$\nu:  \beta \to {\partial f_0\over \partial \beta}|_{f_0^W}.$$
\end{lemma}
\begin{proof}  It is obvious that the two horizontal sequences are identical if $I$ is an isomorphism. 
But the lemma says with the maps $\nu, \mu_3$ defined independently as above, they are still isomorphic.  
Thus it suffices to prove the commutativity for
\begin{equation}\begin{array}{ccc}
 {TX_W(1)\over G(1)}|_{\bar f_0^W}
&\stackrel{\mu_3}\rightarrow & \mathcal O_{\mathbf P^n}(h+1)|_{f_0^W}\\
\downarrow\scriptstyle{I} && \|\\{TW(1)}|_{f_0^W}
&\stackrel{\nu} \rightarrow &\mathcal O_{\mathbf P^n}(h+1)|_{f_0^W}\end{array}\end{equation}
Next we just verify it at each point. This is a tedious but straightforward verification.
Because $f_0$ is smooth, $X_W\to A$ is a smooth map around the point
$z_0=(f_0, x)\in X_W$. Thus $T_{(f_0, x)}X_W\to T_{f_0}A$ is surjective. So we can let
$${\partial \over \partial a_i}-\beta_i$$ be an inverse of ${\partial \over \partial a_i}$ at the point
$z_0=(f_0, x)$.  Let
$$\sigma'=\sum_{i=0}^h x_i({\partial \over \partial a_i}-\beta_i)+y$$
be a representative of an element  in
$${TX(1)\over G(1)}|_{z_0}$$
where the number $x_i$ is the 
coefficient of ${\partial \over \partial a_i}-\beta_i$, and
$y\in (T_{f_0^W/A}(1))|_{z_0}$. \par
By the definition of $\mu_3$, 
\begin{equation} \mu_3(\sigma')=\sum_{i=1}^{h} x_i(L_1\cdots\hat L_i\cdots L_h)|_{z_0}.\end{equation}
If $x$ is not a zero of all $L_i, i\neq 0$, now at this point $z_0$,  $\sigma'$ can be written as 
\begin{eqnarray}  \sigma' =&& {\sum_{i=0}^h x_i L_0\cdots \hat L_i\cdots L_h\over L_1\cdots L_h}|_x\cdot ({\partial \over \partial a_0}-\beta_0)
\\ &&+{1\over L_1\cdots L_h}\sum_{i=1}^{h} x_iL_1\cdots\hat L_i\cdots L_h (L_0\beta_0-L_i\beta_i)|_x
\\&&+y+\sum_{i=1}^h x_i(
{\partial \over \partial a_i}-{L_0\over L_i} {\partial \over \partial a_0})|_x \end{eqnarray}

Modulo $G(1)$, it is just
\begin{eqnarray} {\sigma'\over G(1)}= &&{\sum_{i=0}^h x_i L_0\cdots \hat L_i\cdots L_h\over L_1\cdots L_h}|_x\cdot ({\partial \over \partial a_0}-\beta_0)
\\ &&+{1\over L_1\cdots L_h}\sum_{i=1}^{h} x_iL_1\cdots\hat L_i\cdots L_h (L_0\beta_0-L_i\beta_i)|_x+y|_x \end{eqnarray}
Using the  definition of the isomorphism $I$ (the formula (2.14)),
\begin{eqnarray} I({\sigma'\over G(1)}) = &&{\sum_{i=0}^h x_i L_0\cdots \hat L_i\cdots L_h\over L_1\cdots L_h}|_x
\cdot \beta_0
\\ &&+{1\over L_1\cdots L_h}\sum_{i=1}^{h} x_iL_1\cdots\hat L_i\cdots L_h (L_0\beta_0-L_i\beta_i)|_x+y|_x .\end{eqnarray}
Since the last two terms lie in $ T_{f_0/A} $, we obtain that 
\begin{eqnarray} \nu (I({\sigma'\over G(1)}) &&={\sum_{i=1}^{h} x_iL_1\cdots\hat L_i\cdots L_h\over 
L_1\cdots L_h}|_{x} \cdot {\partial f_0(x)\over \partial \beta_0}\\
&&
( because \quad {\partial F\over \partial a_0}-{\partial F\over \partial \beta_0}=0)\\&&
=\sum_{i=1}^{h} (x_i L_1\cdots\hat L_i\cdots L_h)|_{x}.\end{eqnarray}
Thus 
\begin{equation} \nu (I({\sigma'\over G(1)}))=\mu_3({\sigma'\over G(1)}).\end{equation}
This proves the lemma at this point $z_0$. \par
If $x$ is a zero of $L_i, i\neq 0$, say $L_1(x)=0$. Then
\begin{equation} I({\sigma'\over G(1)})=\sum_{i=0}^h x_i\beta_i +y.\end{equation}

\begin{equation}\nu (I({\sigma'\over G(1)})) = (x_1L_0L_2L_3L_4L_h)|_{x}.\end{equation}
Again $$\nu (I({\sigma'\over G(1)}))=\mu_3({\sigma'\over G(1)}). $$

This proves the lemma. \par
\end{proof}

\section{Deformation of curves to other hypersurfaces}
In this section, we try to use the surjectivity of $P_A^s$  to
obtain some results on ${TX_A(1)\over G(1)}$. We'll denote the pull-back normal bundle over $C$ by
$N_{c_0}V$ for any smooth $V\subset \mathbf P^n$. Also denote  the image of a map
$\mu$ by $Im(\mu)$.

\bigskip

\begin{proposition} Let  $\mathcal L$ be the hyperplane section  bundle
$c_0^\ast(\mathcal O_{\mathbf P^n}(1))$ over $C$. Assume $P_A^s$ is surjective.

\par
(a) If $C$ is a rational curve and
\begin{equation} \{L_i=0\}\cap \{L_i=0\}\cap C_0=\emptyset, i\neq j, \end{equation}
 the map
\begin{equation}\begin{array}{ccc}
\nu^s:  \quad H^0(c_0^\ast(T\mathbf P^n(1)))&\stackrel{\nu^s}\rightarrow H^0(c_0^\ast(\mathcal O_{\mathbf P^n}(h+1)),\end{array}
 \end{equation}
is surjective. 
\par

(b) If \begin{equation} (h+1)h^0(\mathcal L)-h\geq h^0( c_0^\ast(T\mathbf P^n(1)),\end{equation}
and \begin{equation} \{L_i=0\}\cap \{L_i=0\}\cap C_0=\emptyset, i\neq j, \end{equation}
\begin{equation} dim(Im(\nu^s))=(h+1)h^0(\mathcal L)-h.\end{equation} 
 \end{proposition}
\begin{proof} Because the formula (3.1),  the image of $c_0$ completely lies in $W$.
In general we denote the induced morphism on $H^0$ groups by $\phi^s$ if the morphism on the bundles is $\phi$. 
In lemma (2.2), 
pulling back the diagram in the formula (2.18) to $C$, we obtain the commutative diagram

\begin{equation}\begin{array}{ccc}   H^0(\bar c_0^\ast({TX(1)\over G(1)}))&\stackrel{\mu_3^s}\rightarrow
& H^0(c_0^\ast(\mathcal O_{\mathbf P^n}(h+1))) \\
  \downarrow\scriptstyle{I^s} & &\|   \\
 H^0( c_0^\ast(T\mathbf P^n(1)))&\stackrel{\nu^s}\rightarrow & H^0(c_0^\ast(\mathcal O_{\mathbf P^n}(h+1))).\end{array}\end{equation}
where $I^s$ is the isomorphism induced from $I$. 
Thus $$dim(Im(\nu^s))=dim(Im(\mu_3^s)).$$  \par
Adding one more space in the diagram, we obtain
\begin{equation}\begin{array}{ccccc}  H^0(c_0^\ast(TX(1)))&\stackrel{\phi}\rightarrow & H^0(\bar c_0^\ast({TX(1)\over G(1)}))&\stackrel{\mu_3^s}\rightarrow
& H^0(c_0^\ast(\mathcal O_{\mathbf P^n}(h+1))) \\
&& 
  \downarrow\scriptstyle{I^s} & &\|   \\ &&
 H^0( c_0^\ast(T\mathbf P^n(1)))&\stackrel{\nu^s}\rightarrow & H^0(c_0^\ast(\mathcal O_{\mathbf P^n}(h+1))).\end{array}\end{equation}
The sequence in the first row is not exact. 
Using the assumption for the part (a), $C$ is a rational curve,  $H^1(c_0^\ast(G(1)))=0$ (because $c_0^\ast(G(1))$ is a trivial bundle over $\mathbf P^1$).
Thus $\phi$ is surjective. \par
Using the assumption for the part (b), we would like to show $\phi$ is also surjective for
non-zero genus curve. First by the surjectivity of $(P_A^s)$, 
$$c_0^\ast(TX)\simeq \oplus _{h+1} \mathcal O_C\oplus c_0^\ast(Tf_0),$$
where $\mathcal O_C$ is the trivial bundle generated by the
sections $(P_A^s)^{-1}({\partial \over \partial a_i}), i=0, \cdots, h$.
Thus $$H^0(c_0^\ast(TX(1)))\simeq \oplus_{h+1} H^0(\mathcal L)\oplus H^0(c_0^\ast(Tf_0)).$$

By the definition of $\mu_3$, it is easy to see the
image of $\mu_2^s=\mu_3^s\circ \phi$ is just the subspace,
\begin{equation} \{\sum_{i=0}^h x_i  c_0^\ast(L_0\cdots \hat L_i\cdots L_h)\}\subset 
H^0(c_0^\ast(\mathcal O_{\mathbf P^n}(h+1)))\end{equation}
where $x_i\in H^0(c_0^\ast(\mathcal O_{\mathbf P^n}(1)))$ run through all sections. 
Notice that because $$c_0^\ast(L_i), c_0^\ast(L_j), \ for\ i\neq j$$ do not have common zeros, 
$$\sum_{i=0}^h x_i  c_0^\ast(L_0\cdots \hat L_i\cdots L_h)=0$$ if and only if
$$x_i=\epsilon_i L_i, \quad and \ \sum_{i=0}^h\epsilon_i=0,$$
for some complex numbers $\epsilon_i$.
Thus 
$$ dim(Im(\mu_2^s))=(h+1)h^0(\mathcal L)-h. $$
By the assumption for the part (b), 
$$dim(Im(\phi))\geq dim(Im(\mu_2^s))=(h+1)h^0(\mathcal L)-h \geq h^0(c_0^\ast(T\mathbf P^n(1))).$$
Because $I^s$ is an isomorphism, this means that $\phi$ is surjective   \par

Since $\phi$ is surjective, 
$$dim(Im(\nu^s))=(h+1)h^0(\mathcal L)-h.$$ This proves the part (b). If $C$ is a rational curve, it is automatic that 
$$(h+1)h^0(\mathcal L)-h-h^0(\mathcal O_{\mathbf P^n}(h+1))=0.$$ Thus
$\nu^s$ is surjective. So we proved the part (a).
\quad \end{proof}

\bigskip
\begin{proof}{\rm of the theorem  1.1  }: 
Consider the exact sequence
$$\begin{array}{ccccccccc}
0 & \rightarrow & c_0^\ast(Tf_0(1)) &\rightarrow & c_0^\ast(\mathbf P^n(1))
&\rightarrow & c_0^\ast(\mathcal O_{\mathbf P^n}(h+1))& \rightarrow & 0.\end{array}$$
Then we have the long exact sequence
$$\begin{array}{ccccccccc}
H^0(c_0^\ast(\mathbf P^n(1))) & \stackrel{\nu^s}\rightarrow & H^0( c_0^\ast(\mathcal O_{\mathbf P^n}(h+1))) &&&&&& \\
& &\downarrow &&&&&& \\
& &
 H^1(c_0^\ast(Tf_0(1)) &\rightarrow & H^1(c_0^\ast(\mathbf P^n(1)))
& \rightarrow & H^1( c_0^\ast(\mathcal O_{\mathbf P^n}(h+1)))& \rightarrow & 0.\end{array}$$
Thus the $codim(Im(\nu^s))$ is 
\begin{equation} h^1(c_0^\ast(Tf_0(1)))-h^1( c_0^\ast(T\mathbf P^n(1))+h^1(\mathcal L^{h+1}).\end{equation}
Combining the result from  proposition  3.1, we proved the part (b).\par

For the part (a), by the proposition, $\nu^s$ is surjective. Then 
\begin{equation} h^1(c_0^\ast(Tf_0(1)))-h^1( c_0^\ast(T\mathbf P^n(1))+h^1(\mathcal L^{h+1})=0.\end{equation}
Notice  $h^1( c_0^\ast(T\mathbf P^n(1))=h^1(\mathcal L^{h+1})=0$ for a rational curve,  then 
$$ h^1(c_0^\ast(Tf_0(1)))=0.$$
We complete the proof.
\quad\end{proof}

\section{Rational curves in a smooth quintic threefold}  In this section, we apply above theorem  1.1 to rational curves
in  a smooth quintic 3-fold(which is not generic). 
\bigskip

{\bf Example 4.1} Let $n=4, h=5$ and $d=1$ in theorem  1.1. Now $C$ is a line, thus a rational curve. 
Consider the
Fermat quintic $f_0$, and the parameter space $A$ satisfying the condition in formula (1.5).  In this case,  
$$H^1(N_{c_0}f_0(1))=H^0(\mathcal O_{\mathbf P^1}(0)\oplus \mathcal O_{\mathbf P^1}(-4))\neq 0.$$
Our theorem 1.1  says this is an obstruction to the deformations of $C_0$ to other quintics in $A$. Indeed 
 no lines in  the Fermat quintic $f_0$ can deform to a generic quintic by  Albano and Katz's result ([1], Prop. 2.1).  In [1], one can find  the detailed
description of deformations of the pair,  $line\subset Fermat\ quintic$. Our result here is stronger than  Albano and Katz's because
the quintics in $A$ is not generic.
\par
This example shows that our obstruction in the theorem  1.1  is meaningful and non-trivial.

\bigskip

{\bf Example 4.2}  Let $n=4, h=5$ and $g=0$ in theorem 1.1.  Let 
$$f_0=l g_1+q g_2$$
where $l$ is linear and $q$ is quadratic. Assume all $l, q, g_i$ are generic.
Let $C_0$ be a smooth rational curve of degree $d$, lying on the quadratic surface $\{l=q=0\}$. Now $f_0$ is not smooth, but 
there are only 24 singular points. We may assume $f_0$ is smooth along $C_0$. Assume the parameter space 
$A$ satisfies formula (1.5).  Then  the theorem 1.1  should still be valid for such pair $C_0\subset f_0$. 
Apply it the
pair $(C_0, f_0)$. 
$$ H^1(N_{C_0}f_0(1))\simeq H^0(\mathcal O_{\mathbf P^1}(-3d)\oplus \mathcal O_{\mathbf P^1}(d-2)),$$
which is non-zero if $d\neq 1$. Thus if $C_0$ is not a line, then the pair $C_0\subset f_0$ can't deform
to all hypersurfaces in $A$ in the first order. In particular, they are obstructed to deform to all hypersurfaces
in the first order. 

\bigskip

\begin{theorem}  {\em ( Upper-bound of the degrees of summands)}.    Let $f_0\subset \mathbf P^4$ be a smooth quintic threefold. 
 and $C_0\subset \mathbf f_0$ a smooth rational curve. Then
\par
(1) $$ N_{C_0}f_0=\mathcal O_{\mathbf P^1}(k)\oplus \mathcal O_{\mathbf P^1}(-2-k)$$
where $k\geq -1$ is an integer. 
\par
(2) If \begin{equation} \{L_i=0\}\cap \{L_i=0\}\cap C_0=\emptyset, i\neq j, \end{equation}
and $P_A^s$ is surjective,  $$k<d.$$
\end{theorem} 
\begin{proof}  (1).  By the adjunction formula 
$deg( N_{C_0}f_0)=-2$. Since  all vector bundles over $\mathbf P^1$ is decomposable. Then 
the part (1) is proved. \par
(2).    The part (1) says 
\begin{equation} N_{C_0}f_0(1)\simeq \mathcal O_{\mathbf P^1}(k+d)\oplus \mathcal O_{\mathbf P^1}(-2-k+d)\end{equation}
where $k\geq -1$.
Then $ H^1(N_{C_0}f_0(1))\simeq H^0(\mathcal O_{\mathbf P^1}(-k-d-2)\oplus O_{\mathbf P^1}(k-d))=0$.  By theorem  1.1, 
$ H^1(N_{C_0}f_0(1))=0$. Then $k-d<0$, or $k<d$. 
\end{proof}

\bigskip

In this case, we also have the Clemens' conjecture [2] that is equivalent to the assertion that for a generic quintic $f_0$, 
Kodaira's condition $ H^1(N_{C_0}f_0)=0$ (without a twist) is also a necessary condition for $C_0$ to deform to all quintics. Our theorem  1.1  did not prove the Clemens' conjecture because of the twist on
the normal bundle $N_{C_0}f_0$. Instead, we only obtain an upper bound of $k$ above.

\bigskip
\section{Smooth curves in a hypersurface of  a higher dimension} In this section, we apply theorem  1.1 
to hypersurfaces of dimension
$n\geq 3$.

\bigskip

\begin{theorem}  If $C_0$ is a rational curve and $P_E^s$ is surjective at a smooth hypersurface 
$f_0\subset \mathbf P^n$, then
$$h\leq  2n-2,$$
where $h=deg(f_0)$.
\end{theorem}
\medskip

{\bf Remark}   This theorem has an importance in the deformation theory of a pair of varieties.  Even though the  
bound  for $h$ in the theorem is the same as 
Clemens' bound in [4], but we did not assume $f_0$ is generic. 
Therefore our theorem  5.1  is beyond the Clemens' result in [4]. In order to obtain a better bound (better than $2n-2$), the
condition of higher order deformations of the pair must be used. This is indeed the case in [8],  in which Voisin used the integrability of a ``vertical " distribution on the versal subvariety to improve the Clemens' bound  $2n-2$  to $2n-3$. 
So this may  explain, why if the ``Full" 
deformation (or all higher orders) of the pair $C_0, f_0$ exists, the sharp bound is 1 less than $2n-2$, which is proved by Voisin ( [7], [8]). 
 We suspect $h\leq 2n-2$ is the sharp bound under our assumption that $P_E^s$ is surjective. This speculation 
is equivalent to the assertion:\par
There exists a smooth hypersurface $f_0$( not generic) of degree 
$$h=2n-2$$ such that it contains an irreducible rational curve $C_0$ with the surjective $P_E^s$. 

\medskip
Therefore the significance of Clemens' bound $2n-2$ might lie in the existence of the first order deformation of the pair
$C_0\subset f_0$, while
the significance of Voisin's bound $2n-3$ lies in the existence of the ``Full" deformation $C_t\subset f_t$.

\bigskip

\begin{proof} Because $P_E^s$ is surjective,  we choose  generic sections $L_i, i=0, \cdots, h$ 
 in $H^0(\mathcal O_{\mathbf P^n}(1))$ for $A$. Then all conditions in theorem  1.1 are satisfied. 
By theorem 1.1, we obtain that 
$h^1(c_0^\ast(Tf_0(1)))=0$.  By Riemann-Roch, 
\begin{eqnarray}   &&h^1(c_0^\ast(Tf_0(1))) \\
 &&=h^0(c_0^\ast(Tf_0(1)))-\biggl(Ch(c_0^\ast(f_0(1)))\cdot Tod(TC)\biggr)\nonumber\\
 &&=h^0(c_0^\ast(Tf_0(1)))-\biggl(c_1(c_0^\ast(Tf_0(1)))+{n-1\over 2}(TC)\biggr)\nonumber\\
 &&=h^0(c_0^\ast(Tf_0(1)))-\biggl(c_1(c_0^\ast(T\mathbf P^n(1)))-(h+1)d+{n-1\over 2}c_1(TC)\biggr)\nonumber\\
&&=h^0(c_0^\ast(Tf_0(1)))+(h-2n)d+(n-1)(g-1)\nonumber\\
&&=h^0(c_0^\ast(Tf_0(1)))+(h-2n)d-(n-1)=0\nonumber 
\end{eqnarray}

 Since \begin{eqnarray}   && h^0(c_0^\ast(Tf_0(1)))=h^0(N_{c_0}f_0)+h^0(TC(1))\\
&&=h^0(N_{c_0}f_0)+d+3,\nonumber\end{eqnarray}
Formula (5.1) becomes
\begin{equation} (h-2n+1)d+h^0(N_{c_0}f_0))-(n-4)=0. \end{equation}
To show $h\leq 2n-2$, it suffices to prove that 
$$h^0(N_{c_0}f_0))-(n-4)>0.$$

Now we use the assumption $P_E^s$ is surjective. 
Let $q\in C_0$ be a point.  Let $v\in Tf_0|_q$ but not in $TC_0|_q$. 
There is a $GL(n+1)$ action on 
$$\mathbf P(H^0(\mathcal O_{\mathbf P^n}(h)))\times \mathbf P^n,$$
that preserves the universal hypersurface
$$X_E\subset  P(H^0(\mathcal O_{\mathbf P^n}(h)))\times \mathbf P^n.$$
Then we have a submanifold 
$$O=\{(g^{-1} f_0, g(C_0)): g\in GL(n+1)\}\subset X_E.$$
The tangent space $T_{(f_0, q)}O$  of it at $(f_0, q), q\in C_0$ lies in
$$T_{(f_0, q)}X_E\subset  T_{f_0}P(H^0(\mathcal O_{\mathbf P^n}(h)))\times T_q\mathbf P^n .$$
It is clear that $T_{(f_0, q)}O$ contains a vector 
$$(\bar \alpha_\sigma, \sigma)\in T_{f_0}P(H^0(\mathcal O_{\mathbf P^n}(h)))\times T_q\mathbf P^n,$$
such that $\sigma|_q=v$.  We should note $$\alpha_\sigma\in P(H^0(\mathcal O_{\mathbf P^n}(h)))$$
is a hypersurface obtained by apply some action on $f_0$ and $\bar \alpha_\sigma$ represents the directional vector
of the line through $f_0$ and $\alpha_\sigma$. 
Since $\sigma|_q=v\in Tf_0|_q$, 
the hypersurface $\alpha_\sigma$ viewed as polynomial of $\mathbf C^{n+1}$ lies in the maximal ideal 
of the point $q\in \mathbf C^{n+1}$ (view $q$ as a point in $\mathbf C^{n+1}$).  Thus

$$\alpha_\sigma=\sum \overline{x_i Q_i},$$
where $x_i$ is in $H^0(\mathcal O_{\mathbf P^n}(1))$ vanishing at $q$, $Q_i$ is a monomial in  
$H^0(\mathcal O_{\mathbf P^n}(h-1))$ and   $\overline{x_i Q_i}$ denotes the vector in 
$$T_{f_0}P(H^0(\mathcal O_{\mathbf P^n}(h)))$$
that represents the direction of the line in $P(H^0(\mathcal O_{\mathbf P^n}(h)))$ through two points, 
$f_0$ and $x_i Q_i$.  Let $y\in H^0(\mathcal O_{\mathbf P^n}(1))$ such that $y|_q=1$.

$$\sum_i  x_i (\overline {yQ_i} )-y(\overline{x_iQ_i})$$ is in $H^0(TX_E(1)|_{C_0}).$  Hence
$$g_v=\sum_ix_i P_1^s\circ (P_E^s)^{-1}( \overline{yQ_i})-y \sigma$$ must be in 
$$H^0(Tf_0(1)|_{C_0}),$$
and $(g_v)|_q=v$, where $P_1^s$ is the projection map
$$H^0(TX_E|_{C_0})\to H^0( T\mathbf P^n|_{C_0}).$$
 This shows the dimension
of $$H^0(Tf_0(1)|_{C_0}),$$ is at least $n-2$,  because at the point $q$, $\{g_v\}$ span
$$H^0(N_{c_0}f_0(1)))|_q$$ which has dimension $n-2$.  Hence 
$$h^0(N_{c_0}f_0))-(n-4)\geq n-2-(n-4)= 2.$$
We complete the proof.

\par

\end{proof}

\begin{theorem}  Assume all conditions in theorem  1.1., in particular
$P_A^s$ is surjective.  Then either
$$g (h-n+1)\geq (h-2n)d-n+1$$ or

$$g\geq  {d\over 2}+1.$$

\end{theorem} 

{\bf Remark}. If $C_0$ can actually deform to all hypersurfaces ($C_t\subset f_t$ exist for generic hypersurfaces $f_t$) and $h\geq 2n-1$, Clemens has a better bound for $g$ in ([4]),
$$g\geq {1\over 2} (h-2n+1)d+1.$$
But our result is different from Clemens' in many respects.

\bigskip

\begin{proof}  Suppose otherwise, i.e.
$$g (h-n+1)< (h-2n)d-n+1$$ and
$$g<  {d\over 2}+1.$$

By Riemann-Roch and Serre-duality, the inequality $$(h+1)h^0(\mathcal L)-h\geq h^0( c_0^\ast(T\mathbf P^n(1)),$$
is reduced to $g (h-n+1)\leq (h-2n)d-n+1$ which is satisfied. Thus
if $P_A^s$ is surjective,   theorem  1.1 says
$$\sigma(c_0, f_0)=h^1(c_0^\ast(Tf_0(1))-gh=0.$$

Now we calculate $h^1(c_0^\ast(T\mathbf P^n(1))-gh$. Using the condition
$$ h^1(\mathcal L)=h^0(\mathcal L^\ast\otimes K)=0,
 h^1(\mathcal L^2)=h^0((\mathcal L^\ast)^2\otimes K)=0$$
(because $g<  {d\over 2}+1$), 
by the Riemann-Roch,

\begin{eqnarray} && h^1(c_0^\ast(Tf_0(1)))-gh\\
&& =h^0(Tf_0(1))-(Ch(f_0(1))\cdot Tod(TC_0))-hg\nonumber \\
&& > -(Ch(Tf_0(1))\cdot Tod(TC_0))-gh\nonumber\\
&& =-(c_1(c_0^\ast(Tf_0(1)))+{n-1\over 2}(TC_0)-gh\nonumber\\
&& =-(c_1(c_0^\ast(T\mathbf P^n(1))-(h+1)d+{n-1\over 2}c_1(TC_0))-gh\nonumber\\
&&=(h-2n)d+(n-1)(g-1)-gh>0\nonumber\\
\end{eqnarray}
This is a contradition. 
We complete the proof. 
\end{proof}

\section*{Acknowledgments}
We would like to thank H. Clemens for the help and the encouragement, especially for his enlightening 
communication of  theorem  2.1.

\end{document}